\documentclass{rendvuoto}

\usepackage{amsmath}
\usepackage{amsfonts}
\usepackage{amsthm}
\usepackage{amssymb}
\numberwithin{equation}{section}
 \RendicontiPages

\RendicontiHeadPage{1}{GIOVANNI ALESSANDRINI, EDI ROSSET}{SYMMETRY
OF SINGULAR SOLUTIONS}

\summary{We prove radial symmetry of singular solutions to an
overdetermined boundary value problem for a class of degenerate
quasilinear elliptic equations.}

\title{Symmetry of singular solutions of degenerate quasilinear elliptic equations}

\author{Giovanni Alessandrini, Edi Rosset}

\info{Keywords: \textsf{Symmetry, singular solutions, moving planes, $p$-laplacian}
  Work supported in part by MiUR, PRIN no. 2006014115. \\
  Authors' addresses: Giovanni Alessandrini, Dipartimento di Matematica e Informatica,
  Universit\`a degli Studi di Trieste, Italy, e-mail: \textsf{alessang@univ.trieste.it}\\
  Edi Rosset, Dipartimento di Matematica e Informatica,
  Universit\`a degli Studi di Trieste, Italy,
  e-mail: \textsf{rossedi@univ.trieste.it}
  }

\dedicated{In ricordo di Fabio}

\date{xxxxx}

\begin{document}

\maketitle

\section{Introduction} \label{sec:
introduction} We consider solutions $u$ to

\begin{equation}
  \label{eq:equation}
  \textrm{div}(a(|\nabla u|)\nabla u)=0, \qquad\hbox{in } \Omega\setminus
  \{O\},
\end{equation}
which vanish on $\partial \Omega$
\begin{equation}
  \label{eq:Dirichlet}
  u=0, \qquad\hbox{on } \partial\Omega,
\end{equation}
and have a positive singularity at the origin $O$
\begin{equation}
  \label{eq:singularity}
  \lim_{x\rightarrow O}u(x)=M\in(0,+\infty].
\end{equation}
We prove that if $u$ satisfies the overdetermined boundary
condition
\begin{equation}
  \label{eq:Neumann}
  \frac{\partial u}{\partial \nu}=-c, \qquad\hbox{on } \partial\Omega,
\end{equation}
with $c>0$ constant, then $\Omega$ is a ball centered at $O$ and
$u$ is radially symmetric.

To be more specific, we shall assume $\Omega$ to be a bounded
connected open set in $\mathbb{R}^n$, $n\geq 2$, containing  the
origin $O$, and with $C^{2,\alpha}$-smooth boundary
$\partial\Omega$. The nonlinearity $a$ is assumed to be a $C^1$
function {}from $(0,+\infty)$ to $(0,+\infty)$ and satisfying the
degenerate ellipticity condition
\begin{equation}
  \label{eq:ellipticity}
  0<\lambda\leq 1+\frac{sa'(s)}{a(s)}\leq \Lambda, \quad \hbox{for
  every } s>0,
\end{equation}
for some positive constants $\lambda$, $\Lambda$.

Such a class of quasilinear degenerate elliptic equations, which
strictly contains the one of $p$-Laplacian type equations, was
introduced by Lieberman \cite{l:l} and independently, in the two
dimensional case, by Alessandrini, Lupo and Rosset \cite{l:alr}.

With such assumptions the main result of this note is the
following.

\begin{theorem}
  \label{theo:main}
Let $u\in C^{1,\alpha}(\overline{\Omega}\setminus \{O\})$ be a
weak solution to \eqref{eq:equation}, satisfying the conditions
\eqref{eq:Dirichlet}, \eqref{eq:singularity}. If, in addition, $u$
satisfies \eqref{eq:Neumann} then $\Omega$ is a ball centered at
the origin $O$ and $u$ is radially symmetric.
\end{theorem}

Let us observe that in view of the regularity results by Lieberman
\cite{l:l} it is reasonable to treat solutions in the
$C^{1,\alpha}$-class.

Note also that for singular solutions satisfying
\eqref{eq:equation}--\eqref{eq:singularity} the limit $M$ in
\eqref{eq:singularity} may be finite or infinite depending on the
nonlinearity $a$. This fact is particularly evident in the special
case when $a(t)=t^{p-2}$, $p>1$. One readily sees  that when
$p\leq n$ then we have $M=+\infty$, whereas for $p>n$ one must
have $M<+\infty$. See in this respect Kichenassamy and Veron
\cite{l:kv} and, for a detailed study of singular solutions in two
variables we refer to Rosset \cite{l:edi}.

Our proof is based on an adaptation of the well-known method of
moving planes by Alexandrov and Serrin. The adaptation of the
method to degenerate equations was initiated  in \cite{l:ag}. Here
we refer mostly to arguments introduced in \cite{l:a}, however,
the presence of the singularity adds a little further difficulty,
since it may appear, at a first glance, that the method of moving
planes cannot be used after the moving plane has crossed the
singularity. We shall show that if this is the case for a certain
direction $\xi$, then for the opposite direction $-\xi$ the
problem of ``hitting the singularity'' cannot occur.

\begin{remark}We take this opportunity to point out that an \emph{erratum}
is in order in \cite{l:a}. In fact, it is improperly stated there
that the nonlinearity $a$ may depend on $|\nabla u|$ and
\emph{also on $u$}. This is not correct, in fact one should assume
$a=a(|\nabla u|)$ and with this proviso all the statements there
are correct.
\end{remark}

\section{Proof of Theorem \ref{theo:main}} \label{sec:
Proof}

We recall some basic properties of solutions to equation
\eqref{eq:equation} that we shall use repeatedly in our arguments.
Local solutions to the equation
\begin{equation}
  \label{eq:equation_bis}
  \textrm{div}(a(|\nabla u|)\nabla u)=0,
\end{equation}
are obtained as limits in $C^{1,\alpha}$ of solutions $u_\epsilon$
to regularized equations
\begin{equation}
  \label{eq:reg_equation}
\textrm{div}(a_\epsilon(|\nabla u_\epsilon|)\nabla u_\epsilon)=0,
\end{equation}
where $a_\epsilon$ satisfies the same conditions as $a$ and in
addition is $C^\infty$ and $a_\epsilon\geq\epsilon>0$ everywhere.
Consequently $u_\epsilon$ can also be seen as a strong solution to
the non-divergence uniformly elliptic equation
\begin{equation}
  \label{eq:reg_nondiv_equation}
  \sum_{i,j}\left(\delta_{ij}+\frac{a'_\epsilon(|\nabla u_\epsilon|)}
  {|\nabla u_\epsilon|^2a_\epsilon(|\nabla u_\epsilon|)}
  \frac{\partial u_\epsilon}{\partial x_i}
  \frac{\partial u_\epsilon}{\partial x_j}\right) \frac{\partial ^2 u_\epsilon}{\partial x_i
  x_j}=0.
\end{equation}
Consequently solutions to \eqref{eq:equation_bis} inherit some
properties of strong solutions to uniformly elliptic equations.

We quote, in particular, the \emph{Harnack inequality}, that is,
there exists $C=C(\lambda, \Lambda)$ such that: if $u$ solves
\eqref{eq:equation_bis} in $B_R(x_0)$ and $u\geq 0$ then
\begin{equation}
  \label{eq:Harnack}
   \max_{B_{R/2}(x_0)} u\leq C  \min_{B_{R/2}(x_0)} u.
\end{equation}
A further consequence of the use of the regularized solutions is
that solutions to \eqref{eq:equation_bis} satisfy the comparison
principle in the weak form. That is, if $v,u$ solve
\eqref{eq:equation_bis} in a domain $G$ and $v\leq u$ on $\partial
G$ then $v\leq u$ also inside. In addition, if $|\nabla u|+|\nabla
v|>0$ in $\overline{G}$ then the strong version of the comparison
principle holds, that is if $v\leq u$ on $\partial G$ then either
$v\equiv u$ or $v<u$ in $G$.

Let us also observe that,  if $u$ is a solution to
\eqref{eq:equation_bis}, then for every constant $C$, also $C-u$
is a solution. Hence, by the Harnack inequality, one readily
obtains that the solution to
\eqref{eq:equation}--\eqref{eq:singularity} satisfies
\begin{equation}
  \label{eq:0<u<M}
   0<u(x)<M, \quad\hbox{for every } x\in \Omega\setminus\{O\}.
\end{equation}

Let us now introduce the moving plane apparatus. For any direction
$\xi\in \mathbb{R}^n$, $|\xi|=1$, and for any $t\in \mathbb{R}$,
we define the hyperplane
\[
   \Pi_t^\xi=\{x\in \mathbb{R}^n\ |\ x\cdot\xi=t\}.
\]
We denote by $R_t^\xi$ the reflection in $\Pi_t^\xi$, that is
\[
   R_t^\xi x=2(t-x\cdot\xi)\xi+x.
\]
We shall denote
\[
   \left(R_t^\xi u\right)(x)=u\left(R_t^\xi x\right).
\]
If we agree to say that if $x\cdot\xi<t$, $x$ is on the \emph{left
hand side} of $\Pi_t^\xi$, and conversely $x$ is on the
\emph{right hand side} of $\Pi_t^\xi$ if $x\cdot\xi>t$, we denote
by $R_t^\xi\Omega$ the reflection of the part of $\Omega$ which is
on the left hand side of $\Pi_t^\xi$, that is
\[
   R_t^\xi\Omega=\{x\in \mathbb{R}^n\ |\ x\cdot\xi>t, R_t^\xi x \in
   \Omega\}.
\]
Given $\xi$, we fix $\overline{t}$ such that
$R_{\overline{t}}^\xi\Omega=\emptyset$. Letting $t>\overline{t}$
increase, we denote by  $t(\xi)$  the largest number such that
\[
   R_t^\xi\Omega\subset\Omega, \quad \hbox{for every } t\in
   (\overline{t},t(\xi)).
\]
As is well-known since Serrin \cite{l:s}, when $t=t(\xi)$ one of
the following two cases is satisfied
\medskip

  \quad I) $\partial(R_t^\xi\Omega)$ is tangent to $\partial \Omega$ at
  a point $P\not\in \Pi_t^\xi$,
  \quad II) $\partial(R_t^\xi\Omega)$ is tangent to $\partial \Omega$ at
  a point $P\in \Pi_t^\xi$.
\medskip

\noindent Let us consider the family of moving planes associated
to the opposite direction $-\xi$ and the corresponding
reflections. One can easily verify that
\[
   \Pi_t^\xi=\Pi_{-t}^{-\xi}, \quad\hbox{for every } t,
\]
and also
\[
   R_t^\xi=R_{-t}^{-\xi}, \quad\hbox{for every } t.
\]
Now we observe that for every $s<t(\xi)$ we also have
$-s>t(-\xi)$. In fact
\[
R_s^\xi\Omega \varsubsetneqq \Omega\cap\{x\cdot\xi>s\}
\]
and therefore, applying $R_{-s}^{-\xi}$ to both sides,
\[
\Omega\cap\{x\cdot\xi<s \} \varsubsetneqq R_{-s}^{-\xi}\Omega,
\]
that is  $R_{-s}^{-\xi}\Omega$ is \emph{not} contained in $\Omega$
and hence $-s>t(-\xi)$.

Hence, letting $s$ increase to $t(\xi)$, we obtain
\begin{equation}
    \label{eq:t+t}
    t(\xi)+t(-\xi)\leq 0
\end{equation}
Consequently either $t(\xi)=t(-\xi)=0$ or one of the two numbers
$t(\xi)$, $t(-\xi)$ is strictly negative.

If $t(\xi)=t(-\xi)=0$ then, obviously, $\Omega$ is symmetric in
$\Pi_0^\xi=\Pi_0^{-\xi}$.

Assume now $t(\xi)<0$ (the other case $t(-\xi)<0$ being
equivalent). We simplify our notation by posing
\[
\Pi=\Pi_{t(\xi)}^{\xi}, \quad R=R_{t(\xi)}^{\xi}, \quad
G=R_{t(\xi)}^{\xi}\Omega, \quad v=R_{t(\xi)}^{\xi} u.
\]
Since the origin $O$ is on the right hand side of $\Pi$, by
\eqref{eq:0<u<M}, we have that there exists $N$, $0<N<M$ such that
the level set
\[
E=\{x\in\Omega\setminus\{O\}\ |\ u(x)\geq N\}
\]
is strictly on the right hand side of $\Pi$.

Consequently, on $\overline{G}$, $v=R_{t(\xi)}^\xi u<N$. Now we
observe that on $\partial(G\setminus E)$ we have $v\leq u$, in
fact $\partial(G\setminus E)$ can be decomposed as
\[
\partial(G\setminus E)=(\partial G\cap\Pi)\cup (\partial
G\setminus(E\cup\Pi))\cup(\partial E\cap G)
\]
and we have
\[
v=u, \quad \hbox{on }\partial G\cap\Pi,
\]
\[
v=0\leq u, \quad \hbox{on }\partial G\setminus(E\cup\Pi),
\]
\[
v<N\leq u, \quad \hbox{on }\partial E\cap G.
\]
Hence, by the weak comparison principle,
\[
v\leq u, \quad \hbox{in }G\setminus E,
\]
and also, trivially,
\[
v<N\leq u, \quad \hbox{in }G\cap E.
\]
Consequently
\begin{equation}
   \label{eq:v_sotto_u}
v\leq u \quad \hbox{in }G.
\end{equation}
{}From now on we rephrase arguments in \cite{l:a} to prove that
$\Pi$ is a plane of symmetry for $\Omega$.

Let $U$ be an $\epsilon$-neighborhood of $\partial\Omega$ in
$\Omega$, with $\epsilon$ small enough to have $|\nabla u|>0$ in
$U$ and $O\not\in U$. Let $A$ be the connected component of
$(RU)\cap U$ such that $P\in \partial A$. In view of the boundary
conditions \eqref{eq:Dirichlet}, \eqref{eq:Neumann}, we have
\[
(u-v)(P)=0, \quad \nabla(u-v)(P)=0.
\]
Moreover, when case II) occurs, by applying the arguments in
\cite{l:s} to equation \eqref{eq:equation_bis} we also have that
\[
\frac{\partial^2}{\partial \eta^2}(u-v)(P)=0,
\]
for every direction $\eta$. Since $A\subset G\cap U$, $u-v$ is a
non-negative solution of a uniformly elliptic equation in $A$. By
using the Hopf lemma when case I) occurs and its variant due to
Serrin \cite[Lemma 2]{l:s} when case II) occurs and by the strong
comparison principle, it follows that
\begin{equation}
   \label{eq:u=v}
   u=v, \quad\hbox{in } \overline{A}.
\end{equation}
Now, let us prove that $R(\partial\Omega)\subset \partial\Omega$,
obtaining that $\Pi$ is a plane of symmetry of $\Omega$. Assume,
by contradiction, that there exists a point $Q\in
R(\partial\Omega)\setminus\partial\Omega$. Let $\gamma$ be arc in
$R(\partial\Omega)$ joining $Q$ with $P$. One can find a subarc
$\gamma'$ in $\gamma\cap (U\cup \partial\Omega)$ having as
endpoints $P$ and a point $R\in U$. Since in any neighborhood of
$\gamma'$ one can find points of $RU\cap U$, which can be joined
to $P$ through paths inside $RU\cap U$, it follows that
$\gamma'\subset\partial A$. Therefore the set
$R(\partial\Omega)\cap U\cap \partial A$ is non-empty and on such
a set $v=0<u$, contradicting \eqref{eq:u=v}.

We have therefore proved that for any direction $\xi$ there exists
a plane $\Pi$ orthogonal to $\xi$ which is a plane of symmetry for
$\Omega$. Hence $\Omega$ is a ball $B_R(x_0)$. It remains to prove
that $x_0=O$ and that $u$ is radially symmetric.

Let $b=b(s)$ be the inverse function to $ta(t)$, let $K=R^{n-1}c\
a(c)$ and set
\[
v(r)=\int_r^Rb(K\rho^{1-n})d\rho,\quad 0<r\leq R.
\]
One can easily verify  that $v(|x-x_0|)$ solves
\[
  \left\{ \begin{array}{ll}
  \textrm{div} (a(|\nabla
v|)\nabla v)=0, &
  \mathrm{in}\ B_R(x_0)\setminus\{x_0\} ,\\
  &  \\
  v=0, & \mathrm{on}\ \partial B_R(x_0),\\
  &  \\
  \frac{\partial v}{\partial\nu}=-c, & \mathrm{on}\ \partial B_R(x_0).\\
  \end{array}\right.
\]
Now, in $B_R(x_0)\setminus\{O,x_0\}$, $u(x)$ and $v(|x-x_0|)$
satisfy the same elliptic equation and have the same Cauchy data
on $\partial B_R(x_0)$.

We recall that for such an equation the unique continuation
property holds as long as one of the two solutions has
non-vanishing gradient. It follows {}from a standard continuity
argument that $u(x)=v(|x-x_0|)$ for every $x\neq O,x_0$.

Consequently $u$ and $v$ must have the same singular point, that
is $x_0=O$ and finally $u(x)=v(|x|)$. We conclude the proof
observing that the singular value $M$ can be computed as follows

$$M=\int_0^R b\left(\left(\frac{R}{\rho}\right)^{n-1}c\
a(c)\right)d\rho.$$

\nocite{*}
\def\cprime{$'$}
\providecommand{\bysame}{\leavevmode\hbox to3em{\hrulefill}\thinspace}

\end{document}